\documentclass{article}
\usepackage{amsfonts}
\usepackage{amssymb}
\usepackage{amsmath}
\usepackage{srcltx}

\newcommand{\bnom}{\begin{nom}}
\newcommand{\enom}{\end{nom}}
\newcommand{\theo}{{\rm\bf Theorem}}

\newcommand{\defi}{{\rm\bf Definition}}
\newcommand{\ob}{{\rm\bf Remark}}

\newcommand{\corol}{{\rm\bf Corollary}}

\newcommand{\demo}{{\bf Proof}}

\newcommand{\ex}{{\bf Example}}

\newtheorem{nom}{{\!}}[section]
\addtocounter{nom}{0}

\begin{document}

\title{COMMON FIXED POINT THEOREMS FOR PAIRS OF SUBCOMPATIBLE MAPS}
\author{\normalsize H. BOUHADJERA {\footnotesize {\dag}} \&  C.
GODET-THOBIE {\footnotesize {\ddag}} \\
{\normalsize 1-Universit\'{e} europ\'{e}enne de Bretagne, France.}\\
 {\normalsize
2-Universit\'{e} de Bretagne Occidentale,}\\
{\normalsize Laboratoire de Math\'{e}matiques de Brest; Unit\'{e}
CNRS: UMR
6205 }\\
{\normalsize 6, avenue Victor Le Gorgeu, CS 93837, F-29238 BREST
Cedex 3 FRANCE }}
\date{ {\footnotesize E-Mail: \dag hakima.bouhadjera@univ-brest.fr; \ddag
christiane.godet-thobie@univ-brest.fr}}
\maketitle

\begin{abstract}
In this paper, we introduce the new concepts of subcompatibility
and subsequentially continuity which are respectively weaker than
occasionally weakly compatibility and reciprocally continuity.
With them, we establish a common fixed point theorem for four maps
in a metric space which improves a recent result of Jungck and
Rhoades \cite{JuR}. Also we give another common fixed point
theorem for two pairs of subcompatible maps of Gregu\v{s} type
which extends results of the same authors, Djoudi and Nisse
\cite{DjN}, Pathak et al. \cite{PCKM} and others and we end our
work by giving a third result which generalizes results of Mbarki
\cite{Mba} and others.\\


\end{abstract}

\section{Historical introduction and new definitions}

Let $(\mathcal{X},d)$ be a metric space and let $f$ and $g$ be
two maps from $(\mathcal{X},d)$ into itself. $f$ and $g$ are commuting if $fgx=gfx$ for all $x$ in $\mathcal{X}$.

To generalize the notion of commuting maps, Sessa \cite{Ses} introduced the concept of weakly commuting maps. He defines $f$ and $g$ to be weakly commuting if $$d(fgx,gfx)\leq d(fx,gx)$$ for all $x\in \mathcal{X}$. Obviously, commuting maps are weakly commuting
but the converse is not true.

In 1986, Jungck \cite{Jun1} gave more generalized commuting and weakly
commuting maps called compatible maps. $f$ and $g$ above are called
compatible if
$$\lim_{n\rightarrow \infty}d(fgx_{n},gfx_{n})=0\leqno (1)$$ whenever $(x_{n})$ is a sequence in $\mathcal{X}$ such that $\underset{n\rightarrow \infty}\lim fx_{n}=\underset{n\rightarrow \infty}\lim gx_{n}=t$ for some $t\in \mathcal{X}$. Clearly, weakly commuting maps
are compatible, but the implication is not reversible (see \cite{Jun1}).

Afterwards, the same author with Murthy and Cho \cite{JMC} made another
generalization of weakly commuting maps by introducing the concept of
compatible maps of type $(A)$. Previous $f$ and $g$ are said to
be compatible of type $(A)$ if in place of (1) we have two
following conditions: $$\lim_{n\rightarrow \infty}d(fgx_{n},g^{2}x_{n})=0 \text{  and }\lim_{n\rightarrow \infty}d(gfx_{n},f^{2}x_{n})=0.$$

\noindent It is clear to see that weakly commuting maps are compatible of type $(A)$, from \cite{JMC} it follows that the implication is not
reversible.

In their paper \cite{PaK}, Pathak and Khan extended type $(A)$
maps by introducing the concept of compatible maps of type $(B)$ and compared these maps with compatible and
compatible maps of type $(A)$ in normed spaces. To be compatible of type $(B)$, $f$ and $g$ above have to satisfy,
in lieu of condition $(1)$, the inequalities:
\begin{eqnarray*}
\underset{n\rightarrow \infty}{\lim}d(fgx_{n},g^{2}x_{n})
&\leq &\dfrac{1}{2}\left[\underset{n\rightarrow \infty}{\lim}d(fgx_{n},ft)+\underset{n\rightarrow \infty}{\lim}d(ft,f^{2}x_{n})\right] \\
&&\text{and} \\
\underset{n\rightarrow \infty}{\lim}d(gfx_{n},f^{2}x_{n})
&\leq &\dfrac{1}{2}\left[\underset{n\rightarrow \infty}{\lim}d(gfx_{n},gt)+\underset{n\rightarrow \infty}{\lim}d(gt,g^{2}x_{n})\right].
\end{eqnarray*}
It is clear that compatible maps of type $(A)$ are compatible
of type $(B)$, to show that the converse is not true (see \cite{PaK}).

Further, in 1998, Pathak et al. \cite{PCKM} introduced another
generalization of compatibility of type $(A)$ by giving the
concept of compatible maps of type $(C)$. $f$ and $g$ are said
to be compatible of type $(C)$ if they satisfy the two
inequalities:
\begin{eqnarray*}
\lim_{n\rightarrow \infty}d(fgx_{n},g^{2}x_{n}) &\leq &\frac{1}{3}\left[\lim_{n\rightarrow \infty}d(fgx_{n},ft)+\lim_{n\rightarrow \infty}d(ft,f^{2}x_{n})\right. \\
&&\left. +\lim_{n\rightarrow \infty}d(ft,g^{2}x_{n})\right] \\
&&\text{and} \\
\lim_{n\rightarrow \infty}d(gfx_{n},f^{2}x_{n}) &\leq &\frac{1}{3}\left[\lim_{n\rightarrow \infty}d(gfx_{n},gt)+\lim_{n\rightarrow \infty}d(gt,g^{2}x_{n})\right. \\
&&\left. +\lim_{n\rightarrow \infty}d(gt,f^{2}x_{n})\right].
\end{eqnarray*}
The same authors gave some examples to show that compatible maps of type $(C)$ need not be neither compatible nor compatible of type $(A)$ (resp. type $(B)$).

In \cite{PCKL} the concept of compatible maps of type $(P)$ was
introduced and compared with compatible and compatible maps of type $(A)$. $f$ and $g$ are compatible of type $(P)$ if instead
of $(1)$ we have $$\lim_{n\rightarrow \infty}d(f^{2}x_{n},g^{2}x_{n})=0.$$

Note that compatibility, compatibility of type $(A)$ (resp. $(B)$, $(C)$ and $(P)$) are equivalent
if $f$ and $g$ are continuous.

In his paper \cite{Jun2}, Jungck generalized the compatibility, the
compatibility of type $(A)$ (resp. type $(B)$, $(C)$ and $(P)$) by introducing the concept of weak
compatibility. He defines $f$ and $g$ to be weakly compatible if $ft=gt$ for
some $t\in \mathcal{X}$ implies that $fgt=gft$.

It is known that all of the above compatibility notions imply weakly
compatible notion, however, there exist weakly compatible maps which are
neither compatible nor compatible of type $(A)$, $(B)$, $(C)$ and $(P)$ (see \cite{Ali}).

Recently, Al-Thagafi and Shahzad \cite{AlS} weakened the concept of weakly
compatible maps by giving the new concept of occasionally weakly compatible
maps (owc). Two self-maps $f$ and $g$ of a set $\mathcal{X}$ to be owc if and only if
there is a point $x$ in $\mathcal{X}$ which is a coincidence point of $f$
and $g$ at which $f$ and $g$ commute; i.e., there exists a point $x$ in $\mathcal{X}$ such that $fx=gx$ and $fgx=gfx$.

\medskip In this paper, we weaken the above notion by introducing a new concept
called \textbf{subcompatible maps}.

\bnom{\defi} Let $(\mathcal{X},d)$ be a metric space. Maps $f$
and $g:\mathcal{X}\rightarrow \mathcal{X}$ are said to be  {\bf
subcompatible} if and only if there exists a sequence $(x_{n})$ in $ \mathcal{X}$ such that $\underset{n\rightarrow \infty}{\lim}fx_{n}=\underset{n\rightarrow \infty}{\lim}gx_{n}=t$, $t\in \mathcal{X}$ and which satisfy $\underset{n\rightarrow \infty}{\lim}d(fgx_{n},gfx_{n})=0.$
\enom

Obviously, two owc maps are subcompatible, however the converse is not true
in general. The example below shows that there exist subcompatible maps
which are not owc.

\bnom {\ex} {\rm Let $\mathcal{X}=[0,\infty[$ with
the usual metric $d$. Define $f$ and $g$ as follows: $$fx=x^{2} \;
\text{\rm  and } gx=\left\{
\begin{array}{lll}
x+2 &\; \text{\rm if }& x\in [0,4] \cup ]9,\infty[,
\\
x+12 &\; \text{\rm  if }& x\in ]4,9].
\end{array}
\right.
$$}
\enom

\noindent Let $(x_{n})$ be a sequence in
$\mathcal{X}$ defined by $x_{n}=2+\dfrac{1}{n}$ for $n\in \mathbb{N}^{\ast}=\{1,2,\ldots\}$. Then,
\begin{eqnarray*}
\underset{n\rightarrow \infty}\lim fx_{n}=\underset{n\rightarrow
\infty}\lim x_{n}^{2}=4=\underset{n\rightarrow \infty}\lim gx_{n}=\underset{n\rightarrow \infty}\lim (x_{n}+2),
\end{eqnarray*}
and
\begin{eqnarray*}
fgx_{n}=f(x_{n}+2)=(x_{n}+2)^{2}\rightarrow 16
\text{ when }n\rightarrow \infty \\
gfx_{n}=g(x_{n}^{2})=x_{n}^{2}+12\rightarrow 16
\text{ when  } n\rightarrow \infty
\end{eqnarray*}
thus, $\underset{n\rightarrow \infty}\lim d(fgx_{n},gfx_{n})=0$; that is, $f$ and $g$ are subcompatible.

\medskip
\noindent On the other hand, we have $fx=gx$ if and only if $x=2$ and
\begin{eqnarray*}
fg(2)=f(4)=4^{2}=16\, \\
gf(2)=g(4)=4+2=6
\end{eqnarray*}
then, $f(2)=4=g(2)$ but $fg(2)=16\neq 6=gf(2)$, hence maps $f$ and $g$ are not owc. \hfill \mbox{$_{_{\blacksquare}}$}
\medskip

Clearly, we can resume implications between previous notions by
the following list:

\begin{itemize}
\item Commuting maps $\Rightarrow$ \textrm{Weakly commuting maps}

\item Weakly commuting maps $\Rightarrow$ \textrm{Compatible maps}

\item Weakly commuting maps $\Rightarrow $ \textrm{Compatible maps of type $(A)$}

\item Compatible maps of type $(A)$ $\Rightarrow$ \textrm{Compatible
maps of type $(B)$}

\item Compatible maps of type $(A)$ $\Rightarrow$ \textrm{Compatible
maps of type $(C)$}

\item Compatible maps (resp. Compatible of type $(A)$, $(B)$, $(C)$,
$(P)$) $\Rightarrow$ \textrm{Weakly compatible maps}

\item Weakly compatible maps $\Rightarrow $ \textrm{Occasionally
Weakly compatible maps}

\item Occasionally weakly compatible maps $\Rightarrow $ \textbf{Subcompatible maps}.
\end{itemize}

\medskip In his paper \cite{Pan}, Pant introduced the concept of reciprocally continuity as follows: Self-maps $f$ and $g$ of a
metric space $(\mathcal{X},d)$ are reciprocally continuous if and
only if $\underset{n\rightarrow \infty}\lim fgx_{n}=ft$ and
$\underset{n\rightarrow \infty} \lim gfx_{n}=gt$ whenever
$(x_{n}) \subset \mathcal{X}$ is such that $\underset{n\rightarrow
\infty}\lim fx_{n}=\underset{n\rightarrow \infty}\lim
gx_{n}=t\in \mathcal{X}$. Clearly, any continuous pair is
reciprocally continuous but, the converse is not true in general.

\medskip Our second objective here is to introduce a new concept called the
notion of \textbf{subsequentially continuous maps} which weakens the
concepts of continuity and reciprocally continuity given above.

\bnom {\defi} Two self-maps $f$ and $g$ of a metric space $(\mathcal{X},d)$ are said to be {\bf subsequentially continuous} if and only if
there exists a sequence $(x_{n})$ in $\mathcal{X}$ such that $\underset{n\rightarrow \infty}\lim fx_{n}=\underset{n\rightarrow \infty}
\lim gx_{n}=t$ for some $t$ in $\mathcal{X}$ and satisfy $\underset{n\rightarrow \infty}\lim fgx_{n}=ft$ and $\underset{n\rightarrow \infty}\lim gfx_{n}=gt$.
\enom

\noindent If $f$ and $g$ are both continuous or reciprocally
continuous then they are obviously subsequentially continuous. The
next example shows that there exist subsequentially continuous
pairs of maps which are neither continuous nor reciprocally
continuous.

\bnom{\ex} {\rm Let $\mathcal{X}$ be $[0,\infty[$
endowed with the usual metric $d$ and define $f$ and
$g:\mathcal{X}\rightarrow
\mathcal{X}$ by
\begin{equation*}
fx=\left\{
\begin{array}{lll}
1+x &\text{ if }& 0\leq x\leq 1\\
2x-1 &\text{ if } & 1<x<\infty,
\end{array}
\right. \text{\textit{\ }}gx=\left\{
\begin{array}{lll}
1-x &\text{ if }& 0\leq x<1 \\
3x-2 &\text{ if } & 1\leq x<\infty.
\end{array}
\right.
\end{equation*}}
\enom

\noindent Obviously, $f$ and $g$ are discontinuous at $x=1$.

\noindent Let us consider the sequence $x_{n}=\dfrac{1}{n}$ for $n=1,2,\ldots$. We have
\begin{eqnarray*}
\begin{array}{llll} fx_{n}=1+x_{n}\rightarrow 1=t & \rm \text{ when }
n\rightarrow \infty,
 \\
gx_{n}=1-x_{n}\rightarrow 1 & \rm  \text{ when }  n\rightarrow
\infty ,
\end{array}
\end{eqnarray*}
and
\begin{eqnarray*}
fgx_{n}=f(1-x_{n})=2-x_{n}\rightarrow 2=f(1), \\
gfx_{n}=g(1+x_{n})=1+3x_{n}\rightarrow 1=g(1),
\end{eqnarray*}
therefore $f$ and $g$ are subsequentially continuous.

\noindent Now, let $(x_{n})=\left(1+\dfrac{1}{n}\right)$ for $n=1,2,\ldots$. We have
\begin{eqnarray*}
fx_{n}=2x_{n}-1\rightarrow 1=t, \\
gx_{n}=3x_{n}-2\rightarrow 1=t,
\end{eqnarray*}
and
\begin{eqnarray*}
fgx_{n}=f(3x_{n}-2)=6x_{n}-5\rightarrow 1\neq 2=f(1), \\
\end{eqnarray*}
so $f$ and $g$ are not reciprocally continuous. \hfill \mbox{$_{_{\blacksquare}}$}

\medskip Now, we show the interest of these two definitions by
giving three main results.

\section{A general common fixed point theorem}

We begin by a general common fixed point theorem which improves a
result of \cite{JuR}.

\bnom{\theo}\label{sowc1} Let $f$, $g$, $h$ and $k$ be four
self-maps of a metric space $(\mathcal{X},d)$. If pairs of maps $(f,h)$ and $(g,k)$ are subcompatible and reciprocally continuous, then

\noindent (a) $f$ and $h$ have a coincidence point;

\noindent (b) $g$ and $k$ have a coincidence point.

\noindent Further, let $\varphi :(\mathbb{R}^{+})^{6}\rightarrow \mathbb{R}$ be a continuous function satisfying the following condition:

\noindent $(\varphi_{1}):$ $\varphi (u,u,0,0,u,u)>0$ $\forall u>0$.

\noindent We suppose that $(f,h)$ and $(g,k)$ satisfy, for
all $x$ and $y$ in $\mathcal{X}$,

\noindent $(\varphi_{2}):\varphi (d(fx,gy),d(hx,ky),d(fx,hx),d(gy,ky),d(hx,gy),d(ky,fx)) \leq 0$.

\noindent Then, $f$, $g$, $h$ and $k$ have a unique common fixed point.
\enom

\noindent \demo

\noindent Since pairs of maps $(f,h)$ and $(g,k)$ are subcompatible and reciprocally continuous, then, there exist two sequences $(x_{n})$ and $(y_{n})$ in $\mathcal{X}$ such that

\noindent $\underset{n\rightarrow
\infty}\lim fx_{n}=\underset{n\rightarrow \infty}\lim hx_{n}=t$ for
some $t\in \mathcal{X}$ and which satisfy $$\underset{n\rightarrow \infty}
\lim d(fhx_{n},hfx_{n})=d(ft,ht)=0;$$ $\underset{n\rightarrow \infty}
\lim gy_{n}=\underset{n\rightarrow \infty}\lim ky_{n}=z$ for
some $z\in \mathcal{X}$ and which satisfy $$\underset{n\rightarrow
\infty}\lim d(gky_{n},kgy_{n})=d(gz,kz)=0.$$
Therefore $ft=ht$ and $gz=kz$; that is, $t$ is a coincidence point of $f$ and $h$ and $z$ is a coincidence point of $g$ and $k$.

\noindent Now, we prove that $t=z$. Indeed, by inequality $(\varphi_{2})$, we have
\begin{eqnarray*}
&&\varphi (d(fx_{n},gy_{n}),d(hx_{n},ky_{n}),d(fx_{n},hx_{n}), \\
&& d(gy_{n},ky_{n}),d(hx_{n},gy_{n}),d(ky_{n},fx_{n})) \\
&\leq &0.
\end{eqnarray*}
Since $\varphi $ is continuous, taking the limit as $n\rightarrow \infty$ yields
\begin{equation*}
\varphi (d(t,z),d(t,z),0,0,d(t,z),d(z,t)) \leq 0
\end{equation*}
which contradicts $(\varphi_{1})$ if $t\neq z$. Hence $t=z$.

\noindent Also, we claim that $ft=t$. If $ft\neq t$, using $(\varphi_{2})$, we get
\begin{eqnarray*}
&&\varphi (d(ft,gy_{n}),d(ht,ky_{n}),d(ft,ht),\\
&& d(gy_{n},ky_{n}),d(ht,gy_{n}),d(ky_{n},ft)) \\
&\leq &0.
\end{eqnarray*}
Since $\varphi $ is continuous, at infinity, we obtain
\begin{equation*}
\varphi (d(ft,t),d(ft,t),0,0,d(ft,t),d(t,ft)) \leq 0
\end{equation*}
contradicts $(\varphi_{1})$. Hence $t=ft=ht$.

\noindent Again, suppose that $gt\neq t$, using inequality $(\varphi_{2})$, we get
\begin{eqnarray*}
&&\varphi (d(ft,gt),d(ht,kt),d(ft,ht),d(gt,kt),d(ht,gt),d( kt,ft)) \\
&=&\varphi (d(t,gt),d(t,gt),0,0,d(t,gt),d(gt,t)) \leq 0
\end{eqnarray*}
contradicts $(\varphi_{1})$. Thus $t=gt=kt$. Therefore $t=ft=gt=ht=kt$; i.e., $t=z$ is a common fixed point of maps $f$, $g$, $h$ and $k$.

\noindent Finally, suppose that there exists another common fixed point $w$ of maps $f$, $g$, $h$ and $k$ such that $w\neq t$. Then, by inequality $(\varphi_{2})$, we have
\begin{eqnarray*}
&&\varphi (d(ft,gw),d(ht,kw),d(ft,ht),d(gw,kw),d(ht,gw),d(kw,ft)) \\
&=&\varphi (d(t,w),d(t,w),0,0,d(t,w),d(w,t)) \leq 0
\end{eqnarray*}
which contradicts $(\varphi_{1})$. Hence $w=t$. \hfill\mbox{$_{_{\blacksquare}}$}

\medskip If we let in Theorem \ref{sowc1}, $f=g$ and $h=k$, we get the next
corollary:

\bnom{\corol} Let $f$ and $h$ be self-maps of a metric space
$(\mathcal{X},d)$ such that $f$ and $h$ are subcompatible and reciprocally continuous, then, maps $f$ and $h$ have a coincidence point.

\noindent Further let $\varphi:(\mathbb{R}^{+})^{6}\rightarrow
\mathbb{R}$ be a continuous function satisfying condition $(\varphi_{1})$ and
\begin{equation*}
\varphi (d(fx,fy),d(hx,hy),d(fx,hx),d(fy,hy),d(hx,fy),d(hy,fx)) \leq 0
\end{equation*}
for every $x$ and every $y$ in $\mathcal{X}$, then there exists a
unique point $t\in \mathcal{X}$ such that $ft=ht=t$.
\enom

If we put $h=k$, we get the following result:

\bnom{\corol} Let $f$, $g$ and $h$ be three self-maps of
a metric space $(\mathcal{X},d)$. Suppose that pairs of maps $(f,h)$ and $(g,h)$ are subcompatible and reciprocally continuous, then,

\noindent (a) $f$ and $h$ have a coincidence point;

\noindent (b) $g$ and $h$ have a coincidence point.

\noindent Let $\varphi :(\mathbb{R}^{+})^{6}\rightarrow \mathbb{R}$ be a continuous function satisfying condition $(\varphi_{1})$ and
\begin{equation*}
\varphi (d(fx,gy),d(hx,hy),d(fx,hx),d(gy,hy),d(hx,gy),d(hy,fx)) \leq 0
\end{equation*}
for all $x$, $y$ in $\mathcal{X}$, then maps $f$, $g$ and $h$ have a unique common fixed point $t\in \mathcal{X}$.
\enom

Now, with different choices of the real continuous function $\varphi$, we obtain the following corollary which contains several already
published results.

\bnom{\corol} If in the hypotheses of Theorem \ref{sowc1}, we have
instead of $(\varphi_{2})$ one of the following inequalities, for
all $x$ and $y$ in $\mathcal{X}$, then the four maps have a unique
common fixed point
\begin{eqnarray*}
\begin{array} {lll} (a)\text{\textit{\ }}d(fx,gy) & \leq & \alpha \max \{d(hx,ky),d(hx,fx),d(gy,ky), \\
&& \frac{1}{2}(d( hx,gy)+d(ky,fx))\}
\end{array}
\end{eqnarray*} where $\alpha\in ]0,1[$,\\
\begin{eqnarray*}
\begin{array}{ll}
(b)\text{\textit{\ }}d(fx,gy)(1+\alpha d(hx,ky))
 \\
\leq \alpha \max \{d(hx,fx)d(gy,ky),d(hx,gy)d(ky,fx)\} \\
+\beta \max \{d(hx,ky),d(hx,fx),d(gy,ky),\frac{1}{2}(d(hx,gy)+d(ky,fx))\}
\end{array}
\end{eqnarray*}
where $\alpha \geq 0$ and $0<\beta <1$,\\
\begin{equation*}
(c)\text{\textit{\ }}d^{3}(fx,gy) \leq \frac{d^{2}(hx,fx) d^{2}(gy,ky)+d^{2}(hx,gy) d^{2}(ky,fx)}{1+d(hx,ky)+d(hx,fx)+d(gy,ky)},
\end{equation*}
\begin{eqnarray*}
\begin{array}{lll}
(d)\text{\textit{\ }}d(fx,gy) &\leq &\digamma [\max \{d(hx,ky),d(hx,fx), \\
&& d(gy,ky),\frac{1}{2}(d(hx,gy)+d(ky,fx))\}]
\end{array}
\end{eqnarray*}
where $\digamma :\mathbb{R}^{+}\rightarrow \mathbb{R}^{+}$ is an
upper semi-continuous function such that, for every $t>0$,
$0<\digamma (t)<t$.
\enom

\noindent{\demo}

\noindent For proof of $(a)$, $(b)$, $(c)$ and $(d)$, we use
Theorem \ref{sowc1} with
the next functions $\varphi$ which satisfy, for every case, hypothesis $(\varphi_{1})$.

\noindent For $(a)$:
\begin{eqnarray*}
\begin{array}{lll}
&&\varphi (d(fx,gy),d(hx,ky),d(fx,hx),d(gy,ky),d(hx,gy),d(ky,fx)) \\
&=& d(fx,gy)-\alpha\max \{d(hx,ky),d(fx,hx),d(gy,ky), \\
&& \frac{1}{2}(d(hx,gy)+d(ky,fx))\}
\end{array}
\end{eqnarray*}
this function $\varphi$ is used by many authors, for example
Example 1 of Popa \cite{Pop1}.

\noindent For $(b)$:
\begin{eqnarray*}
\begin{array}{lll}
&&\varphi (d(fx,gy),d(hx,ky),d(fx,hx),d(gy,ky),d(hx,gy),d(ky,fx)) \\
&=&(1+\alpha d(hx,ky))d(fx,gy)-\alpha
\max \{d(fx,hx)d(gy,ky), \\
&& d(hx,gy)d(ky,fx)\}-\beta \max \{d(hx,ky),d(fx,hx),d(gy,ky), \\
&& \frac{1}{2}(d(hx,gy)+d(ky,fx))\}
\end{array}
\end{eqnarray*}
for $\beta=1$, we have Example 3 of Popa \cite{Pop2}.

\noindent For $(c)$:
\begin{eqnarray*}
\begin{array}{lll}
&&\varphi (d(fx,gy),d(hx,ky),d(fx,hx),d(gy,ky),d(hx,gy),d(ky,fx)) \\
&=&d^{3}(fx,gy)-\dfrac{d^{2}(fx,hx)d^{2}(gy,ky)+d^{2}(hx,gy)d^{2}(ky,fx)}{1+d(hx,ky)+d(fx,hx)+d(gy,ky)}
\end{array}
\end{eqnarray*}
this function $\varphi$ is the one of Example 5 of \cite{Pop1} with $c=1$.

\noindent For $(d)$:
\begin{eqnarray*}
\begin{array}{lll}
&&\varphi (d(fx,gy),d(hx,ky),d(fx,hx),d(gy,ky),d(hx,gy),d(ky,fx)) \\
&=&d(fx,gy)-\digamma [\max \{d(hx,ky),d(fx,hx),d(gy,ky), \\
&& \frac{1}{2}(d(hx,gy)+d(ky,fx))\}].
\end{array}
\end{eqnarray*}
\hfill \mbox{$_{_{\blacksquare}}$}

\medskip Now, using the recurrence on $n$, we get the following theorem:

\bnom{\theo} Let $h$, $k$ and $\{f_{n}\}_{n\in \mathbb{N}^{\ast}}$ be maps from a metric space $(\mathcal{X},d)$ into itself such that pairs of maps $(f_{n},h)$ and $(f_{n+1},k)$ are subcompatible and reciprocally continuous, then

\noindent (a) $(f_{n},h)$ have a coincidence point;

\noindent (b) $(f_{n+1},k)$ have a coincidence point.

\noindent Suppose that maps $f_{n}$, $f_{n+1}$, $h$ and $k$ satisfy the inequality:
\begin{eqnarray*}
(\varphi_{2})&&\varphi (d(f_{n}x,f_{n+1}y),d(hx,ky),d(f_{n}x,hx), \\
&& d(f_{n+1}y,ky),d(hx,f_{n+1}y),d(ky,f_{n}x)) \\
&\leq &0
\end{eqnarray*}
for all $x$ and $y$ in $\mathcal{X}$, for every $n\in \mathbb{N}^{\ast}$, where $\varphi$ is as in Theorem \ref{sowc1}, then, $h$, $k$ and $\{f_{n}\}_{n\in \mathbb{N}^{\ast}}$ have a unique
common fixed point.
\enom

\noindent{\demo}

\noindent By letting $n=1$, we get the assumptions of Theorem
\ref{sowc1} for maps $h$, $k$, $f_{1}$ and $f_{2}$ with the
unique common fixed point $t$. Now, $t$ is a common fixed point of
$h$, $k$, $f_{1}$ and of $h$, $k$, $f_{2}$. Otherwise, if
$z$ is another common fixed
point of $h$, $k$ and $f_{1}$, then by inequality $(\varphi_{2})$, we have
\begin{equation*}
\left.
\begin{array}{c}
\varphi (d(f_{1}z,f_{2}t),d(hz,kt),d(hz,f_{1}z), \\
d(kt,f_{2}t),d(hz,f_{2}t),d(kt,f_{1}z)) \\
=\varphi (d(z,t),d(z,t),0,0,d(z,t),d(t,z))\leq 0
\end{array}
\right.
\end{equation*}
contradicts $(\varphi_{1})$, then $z=t$.

\noindent By the same manner, we prove that $t$ is the unique common fixed
point of maps $h$, $k$ and $f_{2}$.

\noindent Now, letting $n=2$, we obtain the hypotheses of Theorem
\ref{sowc1} for maps $h$, $k$, $f_{2}$ and $f_{3}$ and then,
they have a unique common fixed point $z$. Analogously,
$z$ is the unique common fixed point of $h$, $k$,
$f_{2}$ and of $h$, $k$, $f_{3}$. Thus $z=t$. Continuing
by this method, we clearly see that $t$ is the required element.
\hfill \mbox{$_{_{\blacksquare}}$}

\bnom{\ob} {\rm We can also have common fixed point by using only
four distances instead of six. The next theorem shows this fact.}
\enom

\bnom{\theo} Let $f$, $g$, $h$ and $k$ be self-maps of a metric
space $(\mathcal{X},d)$. If pairs of maps $(f,h)$ and $(g,k)$ are subcompatible and reciprocally continuous, then,

\noindent (a) $f$ and $h$ have a coincidence point;

\noindent (b) $g$ and $k$ have a coincidence point.

\noindent Let $\psi:(\mathbb{R}^{+})^{4}\rightarrow \mathbb{R}$ be a continuous function such that

\noindent $(\psi_{1}):$ $\psi(u,u,u,u)>0$ $\forall u>0$.

\noindent Suppose that $(f,h)$ and $(g,k)$ satisfy the following
inequality $(\psi_{2})$, for all $x$ and $y$ in
$\mathcal{X}$,

\noindent $(\psi_{2}):$ $\psi(d(fx,gy),d(hx,ky),d(hx,gy),d(ky,fx)) \leq 0$.

\noindent Then, $f$, $g$, $h$ and $k$ have a unique common fixed point.
\enom

\noindent{\demo}

\noindent First, proof of (a) and (b) is similar to proof of first part of Theorem \ref{sowc1}.

\noindent Now, suppose that $d(t,z)>0$, then, using inequality $(\psi_{2})$, we get
\begin{equation*}
\psi(d(fx_{n},gy_{n}),d(hx_{n},ky_{n}),d( hx_{n},gy_{n}),d(ky_{n},fx_{n})) \leq 0.
\end{equation*}

\noindent Since $\psi$ is continuous, we obtain at infinity
\begin{equation*}
\psi(d(t,z),d(t,z),d(t,z),d(z,t)) \leq 0
\end{equation*}
which contradicts $(\psi_{1})$, therefore $z=t$.

\noindent If $d(ft,t)>0$, by inequality $(\psi_{2})$, we have
\begin{equation*}
\psi(d(ft,gy_{n}),d(ht,ky_{n}),d(ht,gy_{n}),d(ky_{n},ft)) \leq 0.
\end{equation*}
Since $\psi$ is continuous, when $n$ tends to infinity, we get
\begin{equation*}
\psi(d(ft,t),d(ft,t),d(ft,t),d(t,ft)) \leq 0
\end{equation*}
which contradicts $(\psi_{1})$, hence $t=ft=ht$.

\noindent Similarly, we have $t=gt=kt$.

\noindent The uniqueness of the common fixed point $t$ follows easily from inequality $(\psi_{2})$ and condition $(\psi_{1})$. \hfill \mbox{$_{_{\blacksquare}}$}

\section{A type Gregu\v{s} common fixed point theorem}

In 1998, Pathak et al. \cite{PCKM} introduced an extension of compatibility
of type $(A)$ by giving the notion of compatibility of type $(C)$ and they proved a common fixed point theorem of Gregu\v{s}
type for four compatible maps of type $(C)$ in a Banach space.
Further, Djoudi and Nisse \cite{DjN} generalized the result of \cite{PCKM}
by weakening compatibility of type $(C)$ to weak compatibility
without continuity. In 2006, Jungck and Rhoades \cite{JuR} extended the
result of Djoudi and Nisse by using an idea called occasional weak
compatibility of Al-Thagafi and Shahzad \cite{AlS} which will be published
in 2008.

\medskip In this part, we establish a common fixed point theorem for four
subcompatible maps of Gregu\v{s} type in a metric space which
extends the results of \cite{DjN}, \cite{JuR} and \cite{PCKM}.

\medskip Let $\mathcal{F}$ be the family of maps $F$ from $\mathbb{R}^{+}$ into itself
such that $F$ is upper semi-continuous and $F(t)<t$ for any $t>0$.

\bnom{\theo}\label{sowc2} Let $f$, $g$, $h$ and $k$ be maps from a metric space $(\mathcal{X},d)$ into itself. If pairs of maps $(f,h)$ and $(g,k)$ are compatible and subsequentially continuous, then,

\noindent (a) $(f,h)$ has a coincidence point;

\noindent (b) $(g,k)$ has a coincidence point.

\noindent Moreover, suppose that the four maps satisfy the following inequality:
\begin{eqnarray*}
\begin{array}{lll}
(2)\text{\textit{\ }}d^{p}(fx,gy) &\leq & F(ad^{p}(hx,ky)+(1-a)\max \{\alpha
d^{p}(fx,hx), \\
&& \beta d^{p}(gy,ky),d^{\frac{p}{2}}(fx,hx)d^{\frac{p}{2}}(fx,ky),\\&& d^{\frac{p}{2}}(fx,ky)d^{\frac{p}{2}}(hx,gy), \\
&& \frac{1}{2}(d^{p}(fx,hx)+d^{p}(gy,ky))\})
\end{array}
\end{eqnarray*}

\noindent for all $x$ and $y$ in $\mathcal{X}$, where $0<a<1$, $\{\alpha ,\beta\} \subset ]0,1]$, $p\in \mathbb{N}^{\ast}$ and $F\in \mathcal{F}$. Then $f$, $g$, $h$ and $k$ have a unique common fixed point.
\end{nom}

\noindent{\demo}

\noindent Since pairs of maps $(f,h)$ and $(g,k)$ are compatible and subsequentially continuous, then, there exist two sequences $(x_{n})$ and $(y_{n})$ in $\mathcal{X}$ such that

\noindent $\underset{n\rightarrow
\infty}\lim fx_{n}=\underset{n\rightarrow \infty}\lim hx_{n}=t$ for
some $t\in \mathcal{X}$ and which satisfy $$\underset{n\rightarrow \infty}
\lim d(fhx_{n},hfx_{n})=d(ft,ht)=0;$$ $\underset{n\rightarrow \infty}
\lim gy_{n}=\underset{n\rightarrow \infty}\lim ky_{n}=z$ for
some $z\in \mathcal{X}$ and which satisfy $$\underset{n\rightarrow
\infty}\lim d(gky_{n},kgy_{n})=d(gz,kz)=0.$$
Therefore $ft=ht$ and $gz=kz$; that is, $t$ is a coincidence point of maps $f$ and $h$ and $z$ is a coincidence point of $g$ and $k$.

\noindent Furthermore, we prove that $t=z$. Suppose that $d(t,z)>0$, indeed by inequality $(2)$ we have
\begin{eqnarray*}
\begin{array}{lll}
d^{p}(fx_{n},gy_{n}) &\leq &F(ad^{p}(hx_{n},ky_{n})+(1-a) \max \{\alpha
d^{p}(fx_{n},hx_{n}), \\
&&\beta d^{p}(gy_{n},ky_{n}),d^{\frac{p}{2}}(fx_{n},hx_{n})d^{\frac{p}{2}}(fx_{n},ky_{n}), \\
&&d^{\frac{p}{2}}(fx_{n},ky_{n})d^{\frac{p}{2}}(hx_{n},gy_{n}), \\
&& \frac{1}{2}(d^{p}(fx_{n},hx_{n})+d^{p}(gy_{n},ky_{n}))\}).
\end{array}
\end{eqnarray*}
By properties of $F$, we get at infinity
\begin{eqnarray*}
d^{p}(t,z) &\leq &F(ad^{p}(t,z)+(1-a)d^{p}(t,z)) \\
&=&F(d^{p}(t,z))<d^{p}(t,z)
\end{eqnarray*}
this contradiction implies that $t=z$.

\noindent Now, if $ft\neq t$, the use of condition $(2)$ gives
\begin{eqnarray*}
\begin{array}{lll}
d^{p}(ft,gy_{n}) &\leq &F(ad^{p}(ht,ky_{n})+(1-a) \max \{\alpha d^{p}(ft,ht), \\
&&\beta d^{p}(gy_{n},ky_{n}),d^{\frac{p}{2}}(ft,ht)d^{\frac{p}{2}}(ft,ky_{n}), \\
&&d^{\frac{p}{2}}(ft,\,ky_{n})d^{\frac{p}{2}}(ht,\,gy_{n}), \\
&&\frac{1}{2}(d^{p}(ft,ht)+d^{p}(gy_{n},ky_{n}))\}).
\end{array}
\end{eqnarray*}
By properties of $F$, we obtain at infinity
\begin{eqnarray*}
\begin{array}{lll}
d^{p}(ft,t) &\leq &F(ad^{p}(ft,t)+(1-a)d^{p}(ft,t)) \\
&=&F(d^{p}(ft,t))<d^{p}(ft,t)
\end{array}
\end{eqnarray*}
this contradiction implies that $t=ft=ht$.

\noindent Similarly, we have $gt=kt=t$. Therefore $t=z$ is a common fixed point of both $f$, $g$, $h$ and $k$.

\noindent Suppose that maps $f$, $g$, $h$ and $k$ have another common fixed point $w\neq t$. Then, by $(2)$ we get
\begin{eqnarray*}
\begin{array}{lll}
d^{p}(ft,gw) &\leq &F(ad^{p}(ht,kw)
+(1-a) \max \{\alpha d^{p}(ft,ht),\beta d^{p}(gw,kw), \\
&&d^{\frac{p}{2}}(ft,ht)d^{\frac{p}{2}}(ft,kw),d^{\frac{p}{2}}(ft,kw)d^{\frac{p}{2}}(ht,gw), \\
&& \frac{1}{2}(d^{p}(ft,ht)+d^{p}(gw,kw))\});
\end{array}
\end{eqnarray*}
that is,
\begin{eqnarray*}
\begin{array}{lll}
d^{p}(t,w) &\leq &F(ad^{p}(t,w)+(1-a) \max \{0,d^{p}(t,w)\}) \\
&=&F(d^{p}(t,w))<d^{p}(t,w)
\end{array}
\end{eqnarray*}
this contradiction implies that $w=t$. \hfill \mbox{$_{_{\blacksquare}}$}

\bnom{\corol} Let $f$ and $h$ be two self-maps of a metric space $(\mathcal{X},d)$. If pair of maps $(f,h)$ is compatible and subsequentially continuous, then $f$ and $h$ have a coincidence point.

\noindent Suppose that $(f,h)$ satisfies the following inequality:
\begin{eqnarray*}
\begin{array}{lll}
d^{p}(fx,fy) &\leq &F(ad^{p}(hx,hy)+(1-a) \max \{\alpha d^{p}(fx,hx),\beta d^{p}(fy,hy), \\
&&d^{\frac{p}{2}}(fx,hx)d^{\frac{p}{2}}(fx,hy),d^{\frac{p}{2}}(fx,hy)d^{\frac{p}{2}}(hx,fy), \\
&& \frac{1}{2}(d^{p}(fx,hx)+d^{p}(fy,hy))\})
\end{array}
\end{eqnarray*}
for all $x$, $y$ in $\mathcal{X}$, where $0<a<1$, $\{\alpha ,\beta\} \subset ]0,1]$, $p\in \mathbb{N}^{\ast}$ and $F\in \mathcal{F}$, then $f$ and $h$ have a unique common fixed point.
\enom

\bnom{\corol} Let $f$, $g$ and $h$ be three self-maps of a metric space $(\mathcal{X},d)$. Suppose that pairs of maps $(f,h)$ and $(g,h)$ are compatible and subsequentially continuous, then,

\noindent (a) $f$ and $h$ have a coincidence point;

\noindent (b) $g$ and $h$ have a coincidence point.

\noindent Further, if the three maps satisfy the next inequality:
\begin{eqnarray*}
\begin{array}{lll}
d^{p}(fx,gy) &\leq &F(ad^{p}(hx,hy)+(1-a) \max \{\alpha d^{p}(fx,hx),\beta d^{p}(gy,hy), \\
&&d^{\frac{p}{2}}(fx,hx)d^{\frac{p}{2}}(fx,hy),d^{\frac{p}{2}}(fx,hy)d^{\frac{p}{2}}(hx,gy), \\
&& \frac{1}{2}(d^{p}(fx,hx)+d^{p}(gy,hy))\})
\end{array}
\end{eqnarray*}
for all $x$ and $y$ in $\mathcal{X}$, where $0<a<1$, $\{\alpha
,\beta \} \subset ]0,1]$, $p\in \mathbb{N}^{\ast}$ and $F\in \mathcal{F}$, then $f$, $g$ and $h$ have a unique common fixed point.
\enom

Again, using the recurrence on $n$, we get the next theorem:

\bnom{\theo} Let $h$, $k$ and $\{f_{n}\}_{n\in \mathbb{N}^{\ast}}$ be self-maps of a metric space $(\mathcal{X},d)$. Suppose that $(f_{n},h)$ and $(f_{n+1},k)$ are compatible and subsequentially continuous, then,

\noindent (a) maps $f_{n}$ and $h$ have a coincidence point;

\noindent (b) $f_{n+1}$ and $k$ have a coincidence point.

\noindent Furthermore, if the maps satisfying the inequality:
\begin{eqnarray*}
\begin{array}{lll}
d^{p}(f_{n}x,f_{n+1}y) &\leq & F \big( ad^{p}(hx,ky)+\\
&&(1-a) \max
\{\alpha d^{p}(f_{n}x,hx),\beta d^{p}(f_{n+1}y,ky),\\
&&d^{\frac{p}{2}}(f_{n}x,hx)d^{\frac{p}{2}}(f_{n}x,ky),d^{\frac{p}{2}}(f_{n}x,ky)d^{\frac{p}{2}}(hx,f_{n+1}y), \\
&& \frac{1}{2}(d^{p}(f_{n}x,hx)+d^{p}(f_{n+1}y,ky))\}
\end{array}
\end{eqnarray*}
for all $x$ and $y$ in $\mathcal{X}$, where $0<a<1$, $\{\alpha
,\beta \} \subset ]0,1]$, $p\in \mathbb{N}^{\ast}$ and $F\in
\mathcal{F}$, then $h$, $k$ and $\{f_{n}\}_{n\in \mathbb{N}^{\ast}}$ have a unique common fixed point.
\enom

\section{A near-contractive common fixed point theorem}

We end our work by establishing the next result which especially
improves the main result of \cite{Mba}.

\bnom{\theo}\label{sowc3} Let $(\mathcal{X},d)$ be a metric space, $f$, $g$, $h$ and $k$ be maps from $\mathcal{X}$ into itself. If $(f,h)$ and $(g,k)$ are compatible and subsequentially continuous or subcompatible and reciprocally continuous, then,

\noindent (a) $f$ and $h$ have a coincidence point;

\noindent (b) $g$ and $k$ have a coincidence point.

\noindent Let $\Phi$ be a continuous function of $[0,\infty[$ into itself such that $\Phi (t)=0$ if and only if $t=0$ and satisfying inequality
\begin{eqnarray*}
( 3) \ \Phi (d(fx,gy))
&\leq &a(d(hx,ky)) \Phi (d(hx,ky)) \\
&&+b(d(hx,ky)) \min \{\Phi (d(hx,gy) ),\Phi (d(ky,fx))\}
\end{eqnarray*}
for all $x$ and $y$ in $\mathcal{X}$, where $a$, $b:[0,\infty[\rightarrow [0,1[$ are upper semi-continuous
and satisfying the condition:
\begin{equation*}
a(t)+b(t)<1 \; \forall t>0.
\end{equation*}
Then $f$, $g$, $h$ and $k$ have a unique common fixed point.
\enom

\noindent{\demo}

\noindent First, proof of parts (a) and (b) is similar to proof of Theorem \ref{sowc1}.

\noindent Now, suppose that $d(t,z)>0$, using inequality $(3)$, we get
\begin{eqnarray*}
&&\Phi (d(fx_{n},gy_{n})) \\
&\leq &a(d(hx_{n},ky_{n}) ) \Phi (d(hx_{n},ky_{n})) \\
&&+b(d(hx_{n},ky_{n})) \min \{\Phi (d(hx_{n},gy_{n})),\Phi (d(ky_{n},fx_{n}))\}.
\end{eqnarray*}
By properties of $\Phi$, $a$ and $b$, we get at infinity
\begin{eqnarray*}
\Phi (d(t,z)) &\leq &[a(d(t,z))+b(d(t,z))] \Phi (d(t,z)) \\
&<&\Phi (d(t,z))
\end{eqnarray*}
which is a contradiction. Hence $\Phi (d(t,z))=0$ which implies that $d(t,z)=0$, thus $t=z$.

\noindent Next, if $ft\neq t$, the use of condition $(3)$ gives
\begin{eqnarray*}
&&\Phi (d(ft,gy_{n})) \\
&\leq &a(d(ht,ky_{n})) \Phi (d(ht,ky_{n})) \\
&&+b(d(ht,ky_{n})) \min \{\Phi (d(ht,gy_{n})),\Phi (d(ky_{n},ft))\}.
\end{eqnarray*}
By properties of $\Phi$, $a$ and $b$, we get at infinity
\begin{eqnarray*}
\Phi (d(ft,t)) &\leq &[a(d(ft,t))+b(d(ft,t))] \Phi (d(ft,t)) \\
&<&\Phi (d(ft,t))
\end{eqnarray*}
this contradiction implies that $\Phi (d(ft,t))=0$
and hence $t=ft=ht$.

\noindent Similarly, we have $gt=kt=t$.

\noindent Now, assume that there exists another common fixed point $w$ of maps $f$, $g$, $h$ and $k$ such that $w\neq t$. By inequality $(3)$ and properties of functions $\Phi$, $a$ and $b$, we obtain
\begin{eqnarray*}
\Phi (d(t,w)) &=&\Phi (d(ft,gw)) \leq a(d(ht,kw)) \Phi (d(ht,kw)) \\
&&+b(d(ht,kw)) \min \{\Phi (d(ht,gw)),\Phi (d(kw,ft))\} \\
&=&[a(d(t,w))+b(d(t,w))] \Phi (d(t,w)) \\
&<&\Phi (d(t,w))
\end{eqnarray*}
this contradiction implies that $\Phi (d(t,w))=0\Leftrightarrow d(t,w)=0$, hence $w=t$.
\hfill \mbox{$_{_{\blacksquare}}$}

\bnom{\ob} {\rm A slight improvement is achieved by replacing inequality $(3)$ in Theorem \ref{sowc3} by the following one:
\begin{eqnarray*}
\Phi (d(fx,gy) ) &\leq &a(d(hx,ky)) \Phi (d(hx,ky)) \\
&&+b(d(hx,ky)) \left[\frac{\Phi ^{\frac{1}{p}}(d(hx,gy))+\Phi ^{\frac{1}{p}}(d(ky,fx))}{2}\right]^{p},
\end{eqnarray*}
with $p$ is an integer such that $p\geq 1$.}
\enom

As particular cases, we immediately obtain the two following
corollaries:

\bnom{\corol} Let $f$ and $h$ be self-maps of a metric space $(\mathcal{X},d)$. Assume that pair of maps $(f,h)$ is
compatible and subsequentially continuous or subcompatible and reciprocally continuous, then, $f$ and $h$ have a coincidence point.

\noindent Further, suppose that pair of maps $(f,h)$ satisfies the inequality:
\begin{eqnarray*}
\Phi (d(fx,fy)) &\leq &a(d(hx,hy)) \Phi (d(hx,hy)) \\
&&+b(d(hx,hy)) \min \{\Phi (d(hx,fy)),\Phi (d(hy,fx))\}
\end{eqnarray*}
for all $x$ and $y$ in $\mathcal{X}$, where $\Phi$, $a$ and $b$ are as in Theorem \ref{sowc3}. Then, $f$ and $h$ have a unique common fixed point.
\enom

\bnom{\corol} Let $f$, $g$, $h:\mathcal{X}
\rightarrow \mathcal{X}$ be maps. If pairs of maps $(f,h)$ and $(g,h)$ are compatible and subsequentially continuous or subcompatible and reciprocally continuous, then,

\noindent (a) $f$ and $h$ have a coincidence point;

\noindent (b) $g$ and $h$ have a coincidence point.

\noindent Moreover, suppose that maps $f$, $g$ and $h$ satisfy the following inequality:
\begin{eqnarray*}
\Phi (d(fx,gy)) &\leq &a(d(hx,hy)) \Phi (d(hx,hy)) \\
&&+b(d(hx,hy) ) \min \{\Phi (d(hx,gy)),\Phi (d(hy,fx))\}
\end{eqnarray*}
for all $x$ and $y$ in $\mathcal{X}$, where $\Phi$, $a$ and $b$ are as in Theorem \ref{sowc3}, then, $f$, $g$ and $h$ have a unique common fixed point.
\enom

We end our work by giving the following result which concern a common fixed
point of a sequence of maps. Its proof is easily obtained from Theorem \ref{sowc3} by recurrence.

\bnom{\theo} Let $(\mathcal{X},d)$ be a metric space, $h$, $k$, $\{f_{n}\}_{n\in \mathbb{N}^{\ast}}$ be maps from $\mathcal{X}$ into itself. If pairs of maps $(f_{n},h)$ and $(f_{n+1},k)$ are
compatible and subsequentially continuous or subcompatible and reciprocally continuous, then,

\noindent (a) $f_{n}$ and $h$ have a coincidence point;

\noindent (b) $f_{n+1}$ and $k$ have a coincidence point.

\noindent Let $\Phi$ be a continuous function of $[0,\infty[$ into itself such that $\Phi (t)=0$ if and only if $t=0$ and satisfying the following inequality:
\begin{eqnarray*}
&&\Phi (d(f_{n}x,f_{n+1}y)) \\
&\leq &a(d(hx,ky)) \Phi (d(hx,ky)) \\
&&+b(d(hx,ky)) \min \{\Phi (d(hx,f_{n+1}y)) ,\Phi (d(ky,f_{n}x))\}
\end{eqnarray*}
for all $x$ and $y$ in $\mathcal{X}$, where $a$ and $b:[0,\infty[\rightarrow [0,1[$ are upper semi-continuous and satisfying the condition
\begin{equation*}
a(t)+b(t)<1 \; \;  \forall t>0.
\end{equation*}
Then, $h$, $k$ and all maps $\{f_{n}\}_{n\in \mathbb{N}^{\ast}}$ have a unique common fixed point.
\enom

\end{document}